\documentclass[11pt, a4paper, leqno]{article}
\usepackage[intlimits]{amsmath}
\usepackage{amsthm,amssymb}
\usepackage{amsfonts}
\usepackage{fancyhdr}
\usepackage{verbatim}
\usepackage{graphicx}
\DeclareMathOperator{\re}{Re}

\DeclareMathOperator{\supp}{supp}
\DeclareMathOperator{\Dom}{Dom}
\newcommand{\bgrad}{\text{\bf grad}\,}

\newcommand{\Acal}{\mathcal{A}}

\newcommand{\Lcal}{\mathcal{L}}
\newcommand{\Mcal}{\mathcal{M}}
\newcommand{\Rbb}[1][]{\mathbb{R}^{#1}}
\newcommand{\Cbb}[1][]{\mathbb{C}^{#1}}
\newcommand{\Nbb}{\mathbb{N}}
\newcommand{\dpar}{\partial}
\newcommand{\drm}{\mathrm{d}}
\newcommand{\erm}{\mathrm{e}}
\newcommand{\muu}{\underline{\mu}{}}
\newcommand{\dsp}{\displaystyle}
\allowdisplaybreaks[4]
\newcommand{\scal}[2]{\ensuremath{\langle #1,#2\rangle}}
\DeclareMathOperator{\spec}{spec}
\DeclareMathOperator{\spp}{\spec_{\text{pt}}}
\DeclareMathOperator{\spe}{\spec_{\text{ess}}}
\DeclareMathOperator{\spd}{\spec_{\text{dis}}}
\numberwithin{equation}{section}
\theoremstyle{plain}
\newtheorem{thm}{Theorem}[section]

\newtheorem{lem}[thm]{Lemma}

\newtheorem{prop}[thm]{Proposition}
\newtheorem{cor}[thm]{Corollary}
\theoremstyle{definition}
\newtheorem{defn}[thm]{Definition}

\theoremstyle{remark}
\newtheorem{rem}[thm]{Remark}

\newcommand{\pd}[2]{\frac{\partial{#1}}{\partial{#2}}}
\newcommand{\bu}{\mathbf{u}}
\newcommand{\bdiv}{\text{\bf div}\,}
\newcommand{\nhat}{\hat{\mathbf{n}}}
\newcommand{\ii}{\rm i}
\newcommand{\const}{\text{const}}
\begin{document}
\title{Existence of eigenvalues of a linear operator pencil in a curved waveguide --- localized shelf
waves on a curved coast}
\author{
E R Johnson\\
\normalsize\small Department of Mathematics, University College London\\
\normalsize\small Gower Street, London WC1E 6BT, U.~K.\\
\normalsize\small email {\sffamily erj@math.ucl.ac.uk}
\and
Michael Levitin\\
\normalsize\small Department of Mathematics, Heriot-Watt University\\
\normalsize\small Riccarton, Edinburgh EH14 4AS, U.~K.\\
\normalsize\small email {\sffamily M.Levitin@ma.hw.ac.uk}
\and
Leonid Parnovski\\
\normalsize\small Department of Mathematics, University College London\\
\normalsize\small Gower Street, London WC1E 6BT, U.~K.\\
\normalsize\small email {\sffamily Leonid@math.ucl.ac.uk}}
%
%
%
\date{%
September 2004
\thanks{The research of M.L. and L.P. was
partially supported by the EPSRC Spectral Theory Network}}
\maketitle
\begin{abstract} \noindent The study of the possibility of existence of the non-propagating, trapped 
continental shelf waves (CSWs)along curved coasts reduces mathematically to a spectral problem for a self-adjoint operator 
pencil in a curved strip. Using the methods developed in the setting 
of the waveguide trapped mode problem, we show that such CSWs exist for a wide class of coast curvature and depth profiles.
\end{abstract}

\

\newpage
\section{Introduction}

Measurements of velocity fields along the coasts of oceans throughout the
world show that much of the fluid energy is contained in motions with periods
of a few days or longer. The comparison of measurements at different places
along the same coast show that in general these low-frequency disturbances
propagate along coasts with shallow water to the right in the northern
hemisphere and to the left in the southern hemisphere. These waves have come
to be known as continental shelf waves (CSWs). The purpose of the present
paper is to demonstrate, using the most straightforward model possible, the
possibility of the non-propagating, trapped CSWs along curved coasts.
The existence of such non-propagating modes would be significant as they
would tend to be forced by atmospheric weather systems, which have similar
periods of a few days, similar horizontal extent, and a reasonably broad
spectrum in space and time. Areas where such modes were trapped would thus
appear to be likely to show higher than normal energy in the low frequency
horizontal velocity field.


    The simplest models for CSWs take the coastal oceans to be inviscid and of
constant density. Both these assumptions might be expected to fail in various regions such as when strong
currents pass sharp capes or when the coastal flow is strongly stratified.
However for small amplitude CSWs in quiescent flow along smooth coasts
viscous separation is negligible.
Similarly most disturbance energy is concentrated in the modes with the least vertical
structure, which are well described by the constant density model \cite{LeBloM78}. The
governing equations are then simply the rotating incompressible Euler
equations. Further, coastal flows are shallow in the sense that the ratio  
of depth to typical horizontal scale is small. Expanding the rotating incompressible Euler equations in
powers of this ratio and retaining only the leading order terms gives the rotating
shallow water equations \cite{Pedlo86}.
\begin{align}
\label{SWE1}
\pd{\bu}{t} + \bu\cdot\bgrad\bu - 2\Omega\mathbf{k}\times\bu &= -g\: \bgrad \widetilde{H}, \\
\pd{\widetilde{H}}{t} + \bdiv [(\widetilde{H}+H)\bu] &= 0.
\label{SWE2}
\end{align}
Here $x$ and $y$ are horizontal coordinates in a frame fixed to the rotating Earth, $\mathbf{k}$ 
is a vertical unit vector, $\bu(x,y,t)$ is the horizontal velocity, $\Omega$ is the (locally constant) vertical 
component of the
Earth's rotation, $g$ is gravitational acceleration, $\widetilde{H}(x,y,t)$ is the vertical
displacement of the free surface and $H(x,y)$ is the local undisturbed fluid depth.

System \eqref{SWE1},\eqref{SWE2} admits waves of two types, denoted Class 1 and Class 2 by \cite{Lamb32}. 
Class 1 are fast high-frequency waves, the rotation modified form of
the usual free surface water waves, although here present only as long,
non-dispersive waves 
with speeds of order $\sqrt{gH}$.
Class 2 waves are slower, low-frequency waves that vanish in the absence of
depth change or in the absence of rotation. It is the Class 2 waves that give
CSWs. They have little signature in the vertical height field $\widetilde{H}(x,y,t)$ and
are observed through their associated horizontal velocity fields \cite{Hamon66}.
The Class 1 waves can be removed from \eqref{SWE1},\eqref{SWE2} by considering the
`rigid-lid' limit, where the external Rossby radius $\sqrt{gH}/2\Omega$ (which gives the relaxation
distance of the free surface) is large compared to the horizontal
scale of the motion. This is perhaps the most accurate of the approximations
noted here, causing the time-dependent term to vanish from \eqref{SWE2} and the
right sides of \eqref{SWE1} to become simple pressure gradients.

For small amplitude waves the nonlinear terms in \eqref{SWE1},\eqref{SWE2} are negligible and cross-differentiating
gives
\begin{align}
\label{lin1}
\pd{\zeta}{t} + 2\Omega\bdiv\bu &= 0, \\
\bdiv (H\bu) &= 0.
\label{lin2}
\end{align}
where $\zeta=\pd{v}{x}-\pd{u}{y}$ is the vertical component of relative vorticity. Equation \eqref{lin2}
is satisfied by introducing the volume flux streamfunction defined though
\begin{equation}
Hu = -\pd{\psi}{y}, \qquad Hv = \pd{\psi}{x},
\end{equation}
allowing \eqref{lin1} to be written as the single equation
\begin{equation}
\bdiv \left( \frac{1}{H}\bgrad\pd{\psi}{t}\right) + 2\Omega\mathbf{k}\cdot\bgrad\psi\times\bgrad(1/H) = 0.
\label{twe}
\end{equation}
Equation \eqref{twe} is generally described as the
topographic Rossby-wave equation or the equation for barotropic CSWs. Many
solutions have been presented for straight coasts, where the coast lies along
$y=0$ (say) and the depth $H$ is a function of $y$ alone (described as rectilinear
topography here) \cite{LeBloM78}. These have shown
excellent agreement with observations of CSWs, as in \cite{Hamon66}. There has
been far less discussion of non-rectilinear geometries, where either the coast or the
depth profile or both are not functions of a single coordinate. Yet interesting
results appear. \cite{StockH86,StockH87} present extensive numerical integrations of a low-order spectral model 
of a rectangular lake with idealized topography. For their chosen depth profiles normal modes can be divided into 
two types: {\em basin-wide modes} which extend throughout the lake and localized {\em bay modes}. These bay modes 
correspond to the high-frequency modes found in a finite-element model of the Lake of Lugarno by \cite{Trosch84} and
observed by \cite{StockHSTZ87}.  \cite{Johns89a,StockJ89,StockJ91} give a variational formulation and describe simplified
quasi-analytical models that admit localized trapped bay modes. However the geometry changes
in these models are large, with the sloping lower boundary terminating abruptly
where it strikes a coastal wall. Further \cite{Johns87a} notes that \eqref{twe} is
invariant under conformal mappings and so any geometry that can be mapped
conformally to a rectilinear shelf cannot support trapped modes. The question thus
arises as to whether it is only for the most extreme topographic changes that
shelf waves can be trapped or whether trapping can occur on smoothly varying
shelves. The purpose of the present paper is to show that this is not so: trapped modes can exist on smoothly curving coasts.

The geometry considered here is that of a shelf of finite width lying along an impermeable coast. Thus sufficiently 
far from the coast the undisturbed fluid depth becomes the constant depth of the open ocean. It is shown in 
\cite{Johns91b} that at the shelf-ocean boundary of finite-width rectilinear shelves the tangential velocity component 
$u$ vanishes for waves sufficiently long compared to the shelf width. The wavelength of long propagating
disturbances is proportional to their frequency which is in turn proportional to the slope of the shelf.
Thus it appears that for sufficiently weakly sloping shelves the tangential velocity component, 
i.e. the normal derivative of the streamfunction, at the shelf-ocean boundary can be made arbitrarily small. 
Here this will be taken as also giving a close approximation to the boundary condition at the shelf-ocean boundary 
when this boundary is no longer straight. The unapproximated boundary condition is that the streamfunction and its 
normal derivative are continuous across the boundary where they match to the decaying solution of Laplace's equation 
(to which \eqref{twe} reduces in regions of constant depth). This gives a linear integral condition along the boundary. 
The unapproximated problem will not be pursued further here. The boundary condition at the coast is simply one of impermeability 
and thus on both rectilinear and curving coasts is simply that the streamfunction vanishes. Now consider flows of the form 
\begin{equation}
\psi(x,y,t)=\re \{ \Phi(x,y) \exp (-2\ii\omega\Omega t\},
\end{equation}
so $\Phi(x,y)$ gives the spatial structure of the flow and $\omega$ its non-dimensional frequency. Then $\Phi$ 
satisfies
\begin{align}
\label{prob1}
\frac{1}{H}\,\Delta\,\Phi + \bgrad\left(\frac{1}{H}\right)\cdot\bgrad\Phi+
\frac{i}{\omega}\mathbf{k}\cdot\left(\bgrad\Phi\times\bgrad\left(\frac{1}{H}\right)\right)=0,& \\
\label{bc1}
\Phi  = 0  \text{ at the coast},& \\
\label{bc2}
\nhat\cdot\bgrad\Phi =0 \text{ at the shelf-ocean boundary},
\end{align}
where $\nhat$ is normal to the shelf-ocean boundary.

Mathematically, we are going to study the existence of trapped modes (i.e., the eigenvalues either
embedded into the essential spectrum or lying in the gap of the essential spectrum) for the problem \eqref{prob1}--\eqref{bc2}
in a curved strip. Similar problems for the Laplace operator have been extensively
studied in the literature --- either in a curved strip, or in a straight strip with an
obstacle, or in a strip with compactly perturbed boundary.
In the case of Laplacians with Dirichlet boundary
conditions these problems are usually called `quantum waveguides'; the Neumann case
is usually referred to as `acoustic waveguides'. The important result concerning
quantum waveguides was established in \cite{ES}, \cite{DE}: in the curved waveguides
there always exist a trapped mode. Later this result was extended to more general settings;
in particular, in \cite{KK} it was shown that in the case of mixed boundary
conditions (i.e., Dirichlet conditions on one side of the strip and Neumann conditions
on the other side) trapped modes exist if the strip is curved `in the direction of the
Dirichlet boundary'.

The case of quantum waveguides is more complicated because any eventual eigenvalues are
embedded into the essential spectrum and are, therefore, highly unstable. Therefore,
it is believed that in general the existence of trapped modes in this case
is due to some sort of the
symmetry of the problem (see \cite{ELV}, \cite{DP}, \cite{APV}).

In the present paper we use the approach similar to the one used in \cite{DE}
and \cite{ELV}; however, we have to modify this approach substantially due to the
fact that we are working with a spectral problem for an operator pencil rather than that for 
an ordinary operator. 

The rest of the paper is organized in the following way. In Section 2, we discuss the rigorous mathematical statement 
of the problem; in Section 3, we study the essential spectrum, and in Section 4, state and prove the main result 
on the existence of a discrete spectrum (Theorem 4.1). In particular we show that a trapped mode always exists of 
all of the following conditions are satisfied: (a) the depth profile $H$ does not depend upon the longitudinal coordinate; 
(b) the channel is curved in the direction of the Dirichlet boundary; (c) the curvature is sufficiently small.

Similar results can be obtained in a straight strip if the depth profile $H$ depends nontrivially upon the longitudinal 
coordinate; we however do not discuss this problem here.

\section{Mathematical statement of the problem}

\subsection{Geometry}

The original geometry is a straight planar strip of width $\delta$:
$$
G_0 = \{(x,y):x\in\Rbb,\ y\in(0,\delta)\}\,.
$$

\begin{figure}[thb!]
\begin{center}
\fbox{\resizebox{0.9\textwidth}{!}{\includegraphics*{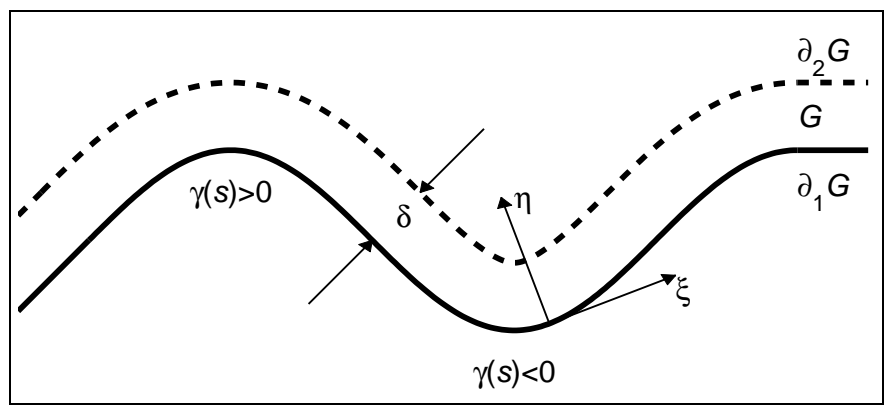}}}
\caption{Domain $G$ and curvilinear coordinates $\xi$, $\eta$. The solid line denotes the boundary 
$\partial_1 G$ with the Dirichlet boundary condition, and the dotted line --- the boundary 
$\partial_2 G$ with the Neumann boundary condition.\label{fig:one}}
\end{center}
\end{figure}

Deformed geometry $G$ is assumed to be a curved planar strip of constant
width $\delta$. To describe it precisely, we introduce the
curve $\Gamma=\{(x=X(\xi),\ y=Y(\xi))\},\,\xi\in\Rbb$, where $\xi$ is a natural arc-length parameter,
i.e., $X'(\xi)^2+Y'(\xi)^2\equiv 1$ . By $\gamma(\xi) = X''(\xi)Y'(\xi) - X'(\xi)Y''(\xi)$ we denote a
(signed) curvature of $\Gamma$ (see Figure~\ref{fig:one}). Note that
$|\gamma(\xi)|^2 = X''(\xi)^2+Y''(\xi)^2$.

We additionally assume
\begin{equation}\label{eq:gam_compact}
\supp\gamma\Subset[-R,R]\qquad\text{for some }R>0\,,
\end{equation}
and set
\begin{equation}\label{eq:gam_ext}
\kappa^+=\sup_{\xi\in [-R,R]}\gamma(\xi)\,,\qquad \kappa^-=-\inf_{\xi\in [-R,R]}\gamma(\xi)\,.
\end{equation}

We shall assume throughout the paper the smoothness condition
\begin{equation}\label{eq:gam_smooth}
\gamma\in C^\infty(\Rbb)\,,
\end{equation}
which can be obviously softened.

Now we can introduce, in a neighbourhood of $\Gamma$,
the curvilinear coordinates $(\xi,\eta)$ as
\begin{equation}\label{eq:curv_coords}
x=X(\xi)-\eta Y'(\xi)\,,\qquad y=Y(\xi)+\eta X'(\xi)\,,
\end{equation}
(where $\eta$ is a distance from a point $(x,y)$ to $\Gamma$)
and describe the \emph{deformed strip} $G$ in these coordinates as
\begin{equation}\label{eq:Omega}
G=G_\gamma=\{(\xi,\eta):\ \xi\in\Rbb,\ \eta\in(0,\delta)\}\,.
\end{equation}

\begin{rem}
As sets of points, $G_\gamma \equiv G_0$ for any $\gamma$, but the metrics are different, see below.
We shall often omit index $\gamma$ if the metric is obvious from the context.
\end{rem}

To avoid self-intersections, we must restrict the width of the strip by natural conditions
\begin{equation}\label{eq:kappa_rest}
\kappa^\pm\le A\delta^{-1}\,,\qquad A = \const \in [0,1)\,.
\end{equation}

Finally, it is an easy computation to show that the Euclidean metric in
the curvilinear coordinates has a form $dx^2+dy^2=gd\xi^2+d\eta^2$, where
$$
g(\xi,\eta) = (1+\eta\gamma(\xi))^2\,.
$$
Later on, we shall widely use the notation
\begin{equation}\label{eq:p}
p(\xi,\eta)=(g(\xi,\eta))^{1/2} = 1+\eta\gamma(\xi)\,.
\end{equation}
Note that in all the volume integrals,
$$
\drm G_\gamma=p(\xi,\eta)\,\drm\xi\,\drm\eta=(1+\eta\gamma(\xi))\,\drm\xi\,\drm\eta=p(\xi,\eta)\,\drm G_0\,.
$$

\subsection{Governing equations}

For a given positive continuously differentiable function $H(\xi,\eta)$ (describing a depth profile), we are looking for a 
function $\Phi(\xi,\eta)$ satisfying \eqref{prob1} with spectral parameter $\omega$.

By substituting
\begin{equation}\label{eq:beta}
\beta(\xi,\eta):=\ln H(\xi,\eta)\,,
\end{equation}
and using explicit expressions for differential operators in curvilinear coordinates, we can re-write
\eqref{prob1} as
\begin{equation}\label{eq:govern_coords}
\begin{split}
&\omega\left(
-\frac{1}{p^2}\,\frac{\dpar^2\Phi}{\dpar\xi^2}
-\frac{\dpar^2\Phi}{\dpar\eta^2}
+
\left(\frac{1}{p^3}\,\frac{\dpar p}{\dpar\xi}+\frac{1}{p^2}\,\frac{\dpar \beta}{\dpar\xi}\right)\,\frac{\dpar \Phi}{\dpar\xi}
+\left(\frac{\dpar \beta}{\dpar\eta}-\frac{1}{p}\,\frac{\dpar p}{\dpar\eta}\right)\,\frac{\dpar \Phi}{\dpar\eta}
\right)
\\
&\qquad =\frac{i}{p}\left(
\frac{\dpar \beta}{\dpar\xi}\,\frac{\dpar \Phi}{\dpar\eta}
-\frac{\dpar \beta}{\dpar\eta}\,\frac{\dpar \Phi}{\dpar\xi}
\right)\,.
\end{split}
\end{equation}

\begin{rem}\label{rem:h}
When deducing equation \eqref{eq:govern_coords}, we have cancelled, on both sides, a common positive
factor $\dsp h(\xi,\eta):=\frac{1}{H(\xi,\eta)}=\erm^{-\beta(\xi,\eta)}$. However, we have to
use this factor when considering corresponding variational equations, in order to keep the
resulting forms symmetric. This leads to a special choice of weighted Hilbert spaces below.
\end{rem}
Further on, we only consider the case of a longitudinally uniform monotone depth profile
\begin{equation}\label{eq:beta_eta}
\beta(\xi,\eta)\equiv \beta(\eta)\,,\qquad\beta'(\eta)>0\,,
\end{equation}
in which case the equation \eqref{eq:govern_coords} simplifies to
\begin{equation}\label{eq:govern_simply}
\begin{split}
&\omega\left(-\frac{1}{p^2}\,\frac{\dpar^2\Phi}{\dpar\xi^2}-\frac{\dpar^2\Phi}{\dpar\eta^2}
+\frac{1}{p^3}\,\frac{\dpar p}{\dpar\xi}\,\frac{\dpar \Phi}{\dpar\xi}
+\left(\beta'-\frac{1}{p}\,\frac{\dpar p}{\dpar\eta}\right)\,\frac{\dpar \Phi}{\dpar\eta}\right) \\
&\qquad =-\frac{i}{p}\,\beta'\,\frac{\dpar \Phi}{\dpar\xi}\,,
\end{split}
\end{equation}
with $\dsp \beta'=\frac{\drm\beta}{\drm\eta}$.

\subsection{Boundary conditions}

Let $\partial_1 G=\{(\xi,0):\xi\in\Rbb\}$ and $\partial_2
G=\{(\xi,\delta):\xi\in\Rbb\}$ denote the lower and the upper
boundary of the strip $G$, respectively. 
Boundary conditions \eqref{bc1}, \eqref{bc2} then become 
\begin{equation}\label{eq:bc}
\left.\Phi\right|_{\dpar_1 G}=\left.\frac{\dpar \Phi}{\dpar\eta}\right|_{\dpar_2 G}=0\,.
\end{equation}
\begin{rem}
If the flow is confined to a channel then the Dirichlet  boundary condition \eqref{bc1} 
applies on both channel walls. This leads to a mathematically different problem which we do not consider in this 
paper.
\end{rem}

\subsection{Function spaces and rigorous operator statement}

We want to discuss the function spaces in which everything acts. Let us denote by
$L_2(G; h)$ the Hilbert space of functions $\phi:G\to\Cbb$ which are square-integrable on $G$
with the weight $\dsp h(\eta)\equiv \frac{1}{H}=\exp(-\beta(\eta))$:
$$
\|\phi\|_{L_2(G;h)}^2=\int_G |\phi(\xi,\eta)|^2\,h(\eta)\,\drm G=
\int_{\Rbb}\int_0^\delta |\phi(\xi,\eta)|^2\,h(\eta)p(\xi,\eta)\,\drm\eta\drm\xi\,<\,\infty\,.
$$
The corresponding inner product will be denoted $\dsp {\scal{\cdot}{\cdot}}_{L_2(G;h)}$.
Similarly we can define the space $L_2(F; h)$ for an arbitrary open subset $F$ of $G$.

Let us formally introduce the operators
$$
\Lcal_\gamma: \Phi\mapsto -\frac{1}{p^2}\,\frac{\dpar^2\Phi}{\dpar\xi^2}-\frac{\dpar^2\Phi}{\dpar\eta^2}
+\frac{1}{p^3}\,\frac{\dpar p}{\dpar\xi}\,\frac{\dpar \Phi}{\dpar\xi}
+\left(\beta'-\frac{1}{p}\,\frac{\dpar p}{\dpar\eta}\right)\,\frac{\dpar \Phi}{\dpar\eta}
$$
and
$$
\Mcal_\gamma: \Phi\mapsto -\frac{i}{p}\,\beta'\,\frac{\dpar \Phi}{\dpar\xi}\,.
$$
(The dependence on $\gamma$ is of course via $p$, see \eqref{eq:p}.)
Then \eqref{eq:govern_simply} can be formally re-written as
\begin{equation}\label{eq:LM}
\omega\Lcal_\gamma \Phi=\Mcal_\gamma\Phi\,,
\end{equation}
or via an {\it operator pencil\/}
\begin{equation}
\Acal_\gamma \equiv \Acal_\gamma(\omega) = \omega\Lcal_\gamma- \Mcal_\gamma
\end{equation}
as
\begin{equation}\label{eq:A}
\Acal_\gamma(\omega) \Phi = 0
\end{equation}
The \emph{domain} of the pencil $\Acal_\gamma$ in the $L_2$-sense is
naturally defined as
\begin{equation}\label{eq:DomA}
\Dom(\Acal_\gamma)=\{\Phi\in H^2(G)\,,\ \Phi\text{ satisfies (\ref{eq:bc})}\}\,,
\end{equation}
where $H^2$ denotes a standard Sobolev space.

On the domain \eqref{eq:DomA}, $\Mcal_\gamma$ is symmetric, and $\Lcal_\gamma$ is symmetric and positive in the sense of 
the scalar product $\dsp {\scal{\cdot}{\cdot}}_{L_2(G;h)}$, with
$$
{\scal{\Lcal_\gamma\Phi}{\Phi}}_{L_2(G;h)}=\int_{\Rbb}\int_0^\delta 
\left(\frac{1}{p}\left|\frac{\dpar\Phi}{\dpar\xi}\right|^2+
p\left|\frac{\dpar\Phi}{\dpar\eta}\right|^2\right)\erm^{-\beta}\,\drm\eta\drm\xi
$$

Later on, we shall use a \emph{weak} (or \emph{variational}) form of \eqref{eq:A}, and shall require
some other
function spaces described below.
Let $F\subseteq G$, and suppose its boundary
is decomposed into two disjoint parts: $\partial F = \partial_1 F \sqcup \partial_2 F$.
We introduce the space
\begin{multline*}
\widetilde{C}_0^\infty(F; \partial_1 F) = \{\phi\in C^\infty(F): \overline{\supp\phi}\cap\partial_1 F =
\emptyset,\\
\text{ and there exists }r>0\text{ such that }
\phi(\xi, \eta)=0\text{ for }(\xi,\eta)\in F,\ |\xi|\ge r\}\,,
\end{multline*}
consisting of smooth functions with compact support vanishing near $\partial_1 F$.
By $\widetilde{H}_0^1(F;\partial_1 F; h)$ we denote its closure with respect to the scalar product
$$
{\scal{\phi}{\psi}}_{\widetilde{H}_0^1(F,\partial_1 F;h)} = {\scal{\phi}{\psi}}_{L_2(F;h)} +
{\scal{\bgrad\phi}{\bgrad\psi}}_{L_2(F;h)}\,.
$$

In what follows we shall study the operators $\Lcal_\gamma$, $\Mcal_\gamma$, and the pencil
$\Acal_\gamma$ from a variational point of view. The details are given in the next Section; here we note only that
from now we understand the expression $\scal{\Lcal_\gamma\Psi}{\Psi}_{L_2(G_\gamma,h)}$ as the quadratic form for the
operator $\Lcal_\gamma$, with the quadratic form domain $\widetilde{H}_0^1(G;\partial_1 G; h)$.

The main purpose of this paper is to study the \emph{spectral
properties} of the operator pencil $\Acal_\gamma$. We recall the
following definitions.

A number $\omega\in\Cbb$ is said to belong to the \emph{spectrum} of $\Acal_\gamma$
(denoted $\spec(\Acal_\gamma)$) if
$\Acal_\gamma(\omega)$ is not invertible.

It is easily seen that in our case the spectrum of $\Acal_\gamma$ is real.

We say that $\omega\in\Rbb$ belongs to the \emph{essential spectrum} of the operator pencil
$\Acal_\gamma$ (denoted $\omega\in\spe(\Acal_\gamma)$) if for this $\omega$ the operator $\Acal_\gamma(\omega)$
is non-Fredholm.

We say that $\omega\in\Cbb$ belongs to the {\it point spectrum\/} of the operator pencil
$\Acal_\gamma$ (denoted $\omega\in\spp(\Acal_\gamma)$),
or, in other words, say that $\omega$ is an {\it eigenvalue\/}, if for this $\omega$ there exists a
non-trivial solution $\Psi\in \Dom(\Acal_\gamma)$ of the problem $\Acal_\gamma(\omega)\Psi=0$.

It is known that the essential spectrum is a closed subset of $\Rbb$, and that any point of the spectrum
outside the essential spectrum is an isolated eigenvalue of finite multiplicity. The set of all such
points is called the \emph{discrete spectrum}, and will be denoted  $\spd(\Acal_\gamma)$
There may, however, exist the points of the spectrum which belong
both to the essential spectrum and the point spectrum.

Our main result (Theorem \ref{thm:main} below) establishes some conditions on the curvature $\gamma$ of the waveguide which
guarantee the existence of eigenvalues of $\Acal_\gamma$.

It is more convenient to deal with problems of this type variationally, and we start the next Section with 
an abstract variational scheme suitable for self-adjoint pencils with non-empty essential spectrum.

\section{Essential spectrum}

\subsection{Variational principle for the essential spectrum}

\begin{defn} We set, for $j\in\Nbb$,
\begin{equation}\label{eq:mu_j}
\mu_{\gamma,j} = \sup_{\substack{U\subset \widetilde{H}_0^1(G;\partial_1 G; h)\\\dim U = j}}
\,\inf_{\substack{\Psi\in U,\,\Psi\ne 0}} 
\frac{\dsp\scal{\Mcal_\gamma\Psi}{\Psi}_{L_2(G_\gamma,h)}}{\dsp\scal{\Lcal_\gamma\Psi}{\Psi}_{L_2(G_\gamma,h)}}
\end{equation}
\end{defn}

As $\dsp\scal{\Lcal_\gamma\Psi}{\Psi}_{L_2(G_\gamma,h)}$ is positive, the right-hand sides of \eqref{eq:mu_j} are well defined,
though the numbers $\mu_{\gamma,j}$ may \emph{a priori} be finite or infinite.

Obviously, for any fixed curvature profile $\gamma$ the numbers $\mu_{\gamma,j}$ form a
non-increasing sequence:
$$
\mu_{\gamma,1}\ge\mu_{\gamma,2}\ge\dots\ge\mu_{\gamma,j}\ge\mu_{\gamma,j+1}\ge\dots\,.
$$

\begin{defn} Denote
\begin{equation}\label{eq:muu}
\muu_\gamma = \lim_{j\to\infty}\mu_{\gamma,j}\,.
\end{equation}
\end{defn}

For general self-adjoint operator pencils the analog of \eqref{eq:muu} may be finite or equal to $\pm\infty$; as we shall 
see below,in our case $\muu_\gamma$ is finite

The following result is a modification, to the case of an abstract self-adjoint linear pencil, of the
general variational principle for a self-adjoint operator with an essential spectrum, see
\cite[Prop. 4.5.2]{EBD}.

\begin{prop}\label{prop:EBD}
Either
\begin{enumerate}
\item[(i)] $\muu_\gamma>-\infty$, and then
$\sup\spe(\Acal_\gamma)=\muu_\gamma$,
\end{enumerate}
or
\begin{enumerate}
\item[(ii)] $\muu_\gamma=-\infty$,
and then $\spe(\Acal_\gamma)=\emptyset$.
\end{enumerate}
Moreover, if $\mu_{\gamma,j}>\muu_\gamma$, then $\mu_{\gamma,j}\in\spd(\Acal_\gamma)$.
\end{prop}

Proposition \ref{prop:EBD} ensures that we can use the variational principle \eqref{eq:mu_j} in order
to find the eigenvalues of the pencil $\Acal_\gamma$ lying \emph{above} the supremum $\muu_\gamma$ of the
essential spectrum.

\subsection{Essential spectrum for the straight strip}

The spectral analysis in the case of a straight strip ($\gamma\equiv 0$) is  rather straightforward
as the problem  admits in this case the separation of variables.

Let us seek the solutions of \eqref{eq:govern_simply},
\eqref{eq:bc} in the case of a straight strip ($\gamma\equiv 0$, and so $p\equiv 1$)
in the form
\begin{equation}\label{eq:phi_def}
\Phi(\xi,\eta)=\phi(\eta)\exp(i\alpha\xi)\,;
\end{equation}
it is sufficient to consider only real values of $\alpha$.

After separation of variables, \eqref{eq:govern_simply}, \eqref{eq:bc} are
written, for each $\alpha$, as a one-dimensional transversal spectral problem
\begin{equation}\label{eq:LM_1d}
\omega(-\phi''+\beta'\phi'+\alpha^2\phi) =\alpha\beta'\phi\,,
\qquad \phi(0)=\phi'(\delta)=0\,,
\end{equation}
Alternatively, introduce operators
$$
\mathfrak{l}_\alpha : \phi\mapsto -\phi''+\beta'\phi'+\alpha^2\phi\,,\qquad
\mathfrak{m}_\alpha : \phi\mapsto\alpha\beta'\phi\,,
$$
and a pencil
$$
\mathfrak{a}_\alpha(\omega) = \omega\mathfrak{l}_\alpha-\mathfrak{m}_\alpha\,,
$$
(again understood in an $L_2((0,\delta);h)$ sense with the domain defined similarly to \eqref{eq:DomA}),
and consider a one-dimensional operator pencil spectral problem $\mathfrak{a}_\alpha(\omega)\phi = 0$.

For a fixed value of $\alpha$, the one-dimensional linear operator
pencil \eqref{eq:LM_1d} has the essential spectrum $\{0\}$ and a discrete spectrum $\spec(\mathfrak{a}_\alpha)$; note
that $\spec(\mathfrak{a}_{-\alpha})=-\spec(\mathfrak{a}_\alpha)$.
Denote, for $\alpha>0$, the top of the spectrum of this transversal problem
by $\omega_\alpha = \sup\spec(\mathfrak{a}_\alpha)$.

\begin{lem}\label{lem:1d_sp} Let $\alpha>0$. Then, under condition \eqref{eq:beta_eta},
\begin{enumerate}
\item[(i)] $\spec(\mathfrak{a}_\alpha)\subset[0, +\infty)$;
\item[(ii)] $0<\omega_\alpha<+\infty$;
\item[(iii)] $\omega_\alpha\to +0$ as $\alpha\to\infty$.
\end{enumerate}
\end{lem}

\begin{proof} By the variational principle analogous to Proposition~\ref{prop:EBD}(i),
$$
\omega_\alpha = \sup_{\substack{\phi\in \widetilde H^1_0((0,\delta),0,h)\\\phi\ne0}} J_\alpha(\phi)\,,
$$
where we set
\begin{equation}\label{eq:vp_1d}
J_\alpha(\phi) =
\frac{\dsp\scal{\mathfrak{m}_\alpha\phi}{\phi}_{L_2((0,\delta),h)}}
{\dsp\scal{\mathfrak{l}_\alpha\phi}{\phi}_{L_2((0,\delta),h)}}
= \frac
{\dsp \int_0^\delta \alpha \beta'(\eta)\phi(\eta)^2
h(\eta)\,\drm\eta}{\dsp \int_0^\delta
(-\phi''(\eta)+\beta'(\eta)\phi'(\eta)+\alpha^2\phi(\eta))\phi(\eta)h(\eta)\,\drm\eta}
\end{equation}

After integrating by parts using $\dsp h(\eta)=\erm^{-\beta(\eta)}$, and inverting the quotient,
we get
\begin{equation}\label{eq:Jalp}
J_\alpha(\phi) = \left(\alpha J^{(1)}(\phi)+\frac{1}{\alpha} J^{(2)}(\phi)\right)^{-1}\,,
\end{equation}
where we denote
$$
 J^{(1)}(\phi) = \frac{\dsp\int_0^\delta \erm^{-\beta(\eta)}|\phi(\eta)|^2\,\drm\eta}{\dsp\int_0^\delta
\beta'(\eta)\erm^{-\beta(\eta)}|\phi(\eta)|^2\,\drm\eta}
$$
and
$$
 J^{(2)}(\phi) = \frac{\dsp\int_0^\delta\erm^{-\beta(\eta)}|\phi'(\eta)|^2\,\drm\eta}{\dsp\int_0^\delta
\beta'(\eta)\erm^{-\beta(\eta)}|\phi(\eta)|^2\,\drm\eta}
$$
The statements (ii) and (iii) of the Lemma now follow immediately from the estimates
$$
J^{(1)}(\phi)\ge
\frac{\dsp \inf_{\eta\in(0,\delta)}\erm^{-\beta(\eta)}}{\dsp \sup_{\eta\in(0,\delta)}(\beta'(\eta)\erm^{-\beta(\eta)})}
$$
and
$$
J^{(2)}(\phi)\ge\frac{\pi^2}{4\delta^2}\,
\frac{\dsp \inf_{\eta\in(0,\delta)}\erm^{-\beta(\eta)}}{\dsp
\sup_{\eta\in(0,\delta)}(\beta'(\eta)\erm^{-\beta(\eta)})}\,,
$$
where the latter inequality uses the variational principle and the fact that the principle eigenvalue of the
mixed Dirichlet--Neumann problem spectral problem for the operator $-\dsp\frac{\drm^2}{\drm \eta^2}$
on the interval $(0,\delta)$ is equal to $\dsp \frac{\pi^2}{4\delta^2}$. The statement (i) follows 
from the positivity of the right-hand side of \eqref{eq:Jalp}.
\end{proof}

We are now able to find the essential spectrum of the problem in a
straight strip.

\begin{lem}\label{lem:sp_str}
Assume that conditions \eqref{eq:beta_eta} hold. Then
\begin{equation}\label{eq:speB}
\spe(\Acal_0)=[-\Omega_*,\Omega_*]\,,
\end{equation}
where
\begin{equation}\label{eq:Omstar}
\Omega_* = \sup_{\phi\in \widetilde H^1_0((0,\delta),0,h)}
\frac{\dsp\frac{1}{2}\int_0^\delta
\beta'(\eta)\erm^{-\beta(\eta)}|\phi(\eta)|^2\,\drm\eta}
{\sqrt{\dsp\int_0^\delta\erm^{-\beta(\eta)}|\phi'(\eta)|^2\,\drm\eta\cdot
\dsp\int_0^\delta \erm^{-\beta(\eta)}|\phi(\eta)|^2\,\drm\eta}}>0\,.
\end{equation}
\end{lem}

\begin{proof} It is standard that
$$
\spe(\Acal_0)=\overline{\bigcup_{\alpha\in\Rbb} \spec(\mathfrak{a}_\alpha)}\,.
$$
Thus, by Lemma \ref{lem:1d_sp}, and with account of the
anti-symmetry of the spectrum of $\mathfrak{a}_\alpha$ with
respect to $\alpha$ and its positivity for $\alpha>0$, we have
$$
\sup\spe(\Acal_0)=\sup_{\alpha>0}\omega_\alpha=
\sup_{\alpha>0}\sup_{\phi\in \widetilde H^1_0((0,\delta),0,h)} J_\alpha(\phi)\,.
$$

By maximizing first with respect to $\alpha$, we obtain, from \eqref{eq:Jalp},
$$
J_\alpha(\phi)\le J_{\alpha_*(\phi)}(\phi)
$$
with the maximizer
$$
\alpha_*(\phi) = \sqrt{\frac{J^{(2)}(\phi)}{J^{(1)}(\phi)}}\,.
$$

Maximizing now with respect to $\phi$ gives
$\sup\spe(\Acal_0)=\Omega_*$, with $\Omega_*$ given by \eqref{eq:Omstar}.
Finally, the statement follows from the fact that $\omega_\alpha$ depends continuously on $\alpha>0$ and
thus takes all the values in $(0, \Omega_*]$ by Lemma~\ref{lem:1d_sp}.
\end{proof}

\subsection{Essential spectrum for a curved strip}

It is now a standard procedure to show that under our conditions the essential spectrum of the problem in a
curved strip coincides with the essential spectrum of the problem in a straight strip given by
Lemma~\ref{lem:sp_str}. Namely, we have the following

\begin{lem}\label{lem:sp_cur}
Let us assume the conditions \eqref{eq:gam_compact}, \eqref{eq:gam_smooth} and \eqref{eq:beta_eta} hold.
Then
$$
\spe(\Acal_\gamma)=\spe(\Acal_0)=[-\Omega_*,\Omega_*]\,
$$
with $\Omega_*$ given by \eqref{eq:Omstar}.
\end{lem}

The proof is based on the fact that any solution of the problem \eqref{eq:govern_simply}, \eqref{eq:bc} with $\gamma\not\equiv 0$
(and thus $p\not\equiv 1$) should coincide in
$G\cap\{|\xi|>R>\max(|\inf\supp\gamma|,|\sup\supp\gamma|)\}$ with a
solution of the
same problem for $\gamma\equiv0$.
An analogous result has been proved in a number of similar situations elsewhere, see e.g. \cite{ES, ELV, DP}, so 
we omit the details of the proof.
We briefly note that the inclusion $\spe(\Acal_\gamma)\subseteq\spe(\Acal_0)$ is proved using the separation of variables
as above and a construction of appropriate Weyl's sequences, and in order to prove the inclusion
$\spe(\Acal_\gamma)\supseteq\spe(\Acal_0)$ one can use the Dirichlet-Neumann bracketing and the discreteness of the spectrum of the problem
\eqref{eq:govern_simply}, \eqref{eq:bc} considered in $G\cap\{|\xi|<R\}$ with additional Dirichlet or Neumann boundary conditions imposed on the ``cuts'' $\{\xi=\pm R\}$.

\section{Main result}

Our main result consists in stating some sufficient conditions on the depth profile $\beta(\eta)$ and
the curvature profile $\gamma(\xi)$ which guarantee the existence of an eigenvalue of the pencil $\Acal_\gamma$ lying
outside the essential spectrum.

\begin{thm}\label{thm:main}
Assume, as before, that the condition \eqref{eq:beta_eta} holds.
Assume additionally that
\begin{equation}\label{eq:beta_2diff}
\beta''(\eta)<0\qquad\text{for }\eta\in(0,\delta)\,.
\end{equation}
Then there exists a constant $C_\beta>0$, which depends only on the depth profile $\beta$, such that
$\spd(\Acal_\gamma)\ne\emptyset$ whenever $\gamma$ satisfies conditions
\eqref{eq:gam_compact}, \eqref{eq:gam_smooth} and
\begin{equation}\label{eq:main_cond}
\int \gamma(\xi)\drm\xi > C_\beta \int \gamma(\xi)^2\drm\xi
\end{equation}
\end{thm}

We give an explicit expression for $C_\beta$ below, see \eqref{eq:C_beta}.

An integral sufficient condition \eqref{eq:main_cond} may be
replaced by a pointwise, although  more restrictive, condition:

\begin{cor}\label{cor:main}
Assume that the conditions \eqref{eq:beta_eta} and \eqref{eq:beta_2diff}
hold. Then there exists a constant $\dsp c_{\beta, R}=\frac{C_\beta}{2R}$ which depends only on the 
depth profile $\beta$ and a given $R>0$
such that $\spd(\Acal_\gamma)\ne\emptyset$ whenever $\gamma\not\equiv 0$ satisfies conditions
\eqref{eq:gam_compact}, \eqref{eq:gam_smooth} and
\begin{equation}\label{eq:cor_cond}
0\le\gamma(\xi)<c_{\beta, R}\qquad\text{ for }|\xi|\le R\,.
\end{equation}
\end{cor}

We prove Theorem~\ref{thm:main} using a number of simple lemmas, the central of which
is the following
\begin{lem}\label{lem:test_fn}
Suppose there exists a function $\widetilde\Psi\in \widetilde H^1_0(G_\gamma,h)$ such that
\begin{equation}\label{eq:main_ineq}
\frac{\dsp\scal{\Mcal_\gamma\widetilde\Psi}{\widetilde\Psi}_{L_2(G_\gamma,h)}}
{\dsp\scal{\Lcal_\gamma\widetilde\Psi}{\widetilde\Psi}_{L_2(G_\gamma,h)}} > \Omega_*\,.
\end{equation}
Then there exists $\omega>\Omega_*$ which belongs to $\spd(\Acal_\gamma)$.
\end{lem}

Lemma \ref{lem:test_fn} is just a re-statement of the variational principle of Proposition~\ref{prop:EBD}.
The main difficulty in its application is of course the choice of an appropriate test function $\widetilde\Psi$.
However such choice becomes much easier if we use the following modification of this Lemma which allows us
to consider test functions which are not necessarily square integrable on $G_\gamma$.

Denote, for brevity, $G^r_\gamma=G_\gamma\cap\{|\xi|<r\}$.

\begin{lem}\label{lem:test_fn_mod}
Suppose there exists a function $\Psi$ and a constant $D$ such that, for any $r>R$, we have
$\Psi \in \widetilde H^1_0(G^r_\gamma,h)$ and
\begin{equation}\label{eq:main_ineq_mod}
\scal{\Mcal_\gamma\Psi}{\Psi}_{L_2(G^r_\gamma,h)}- \Omega_*\scal{\Lcal_\gamma\Psi}{\Psi}_{L_2(G^r_\gamma,h)} \ge D > 0\,.
\end{equation}
Then there exists $\omega>\Omega_*$ which belongs to $\spd(\Acal_\gamma)$.
\end{lem}

The proof of Lemma~\ref{lem:test_fn_mod} uses the construction of an appropriate cut-off function $\chi(\xi)$
such that $\widetilde\Psi=\chi\Psi$ satisfies the conditions of Lemma~\ref{lem:test_fn}, cf. \cite[Prop.~1]{DP}.

We now proceed as follows.

Let $\phi_*(\eta)$ be a minimizer in \eqref{eq:Omstar}, and set
$$
\Psi(\xi,\eta) = \phi_*(\eta) \erm^{i\alpha_{\bullet} \xi}\,,
$$
where
\begin{equation}\label{eq:alpha_crit}
\alpha_{\bullet}=\alpha_*(\phi_*) = \sqrt{\frac{J^{(2)}(\phi_*)}{J^{(1)}(\phi_*)}}\,.
\end{equation}

It is important to note that $\Psi$ is in fact an ``eigenfunction'' of the essential spectrum
of $\Acal_\gamma$ corresponding to its highest positive point $\Omega_*$ and that $\phi_*$ is an
eigenfunction of \eqref{eq:LM_1d} with  $\alpha=\alpha_{\bullet}$
(i.e. of the pencil $\mathfrak{a}_{\alpha_{\bullet}}$) again corresponding to the eigenvalue  $\Omega_*$, and so
\begin{equation}\label{eq:phi_star_eq}
\phi''_* = \beta'\phi'_*+(\alpha_{\bullet}^2-\Lambda_*\alpha_{\bullet}\beta')\phi_*\,,\qquad
\phi_*(0)=\phi'_*(\delta)=0
\end{equation}
with $\dsp \Lambda_*:=\frac{1}{\Omega_*}$ 
(cf. \eqref{eq:LM_1d}).

For future use, we summarize the relations obtained so far:
\begin{align*}
\Lcal_0\Psi &=  (-\phi_*''(\eta)+\beta'(\eta)\phi_*'(\eta)+\alpha_{\bullet}\phi_*(\eta))\,\erm^{i\alpha_{\bullet} \xi} =
(\mathfrak{l}_{\alpha_{\bullet}}\phi_*)\,\erm^{i\alpha_{\bullet} \xi}\,,
\\
\Mcal_0\Psi &=  \alpha_{\bullet}\beta'(\eta)\phi_*(\eta)\,\erm^{i\alpha_{\bullet} \xi}
= (\mathfrak{m}_{\alpha_{\bullet}}\phi_*)\,\erm^{i\alpha_{\bullet} \xi}\,,
\\
\Lcal_\gamma\Psi &=  \left(-\phi_*''(\eta)
+\left(\beta'(\eta)-\frac{1}{p(\xi,\eta)}\,\frac{\dpar p(\xi,\eta)}{\dpar\eta}\right)\phi_*'(\eta)
\right.
\\
&\qquad\qquad\left.+\left(\frac{i\alpha_{\bullet}}{p(\xi,\eta)^3}\,\frac{\dpar p(\xi,\eta)}{\dpar\xi}+
\frac{\alpha_{\bullet}^2}{p(\xi,\eta)^2}\right)\,\phi_*(\eta)
\right)\,\erm^{i\alpha_{\bullet} \xi} \,,
\\
\Mcal_\gamma\Psi &=  \frac{\alpha_{\bullet}}{p(\xi,\eta)}\beta'(\eta)\phi_*(\eta)\,\erm^{i\alpha_{\bullet} \xi} \,,
\\
p(\xi,\eta)&=1+\gamma(\xi)\eta\,,
\end{align*}
(with ${}'$ denoting differentiation with respect to $\eta$).

It is important to note that for any $r>0$
$$
\frac{\dsp\scal{\Mcal_0\Psi}{\Psi}_{L_2(G^r_0,h)}}{\dsp\scal{\Lcal_0\Psi}{\Psi}_{L_2(G^r_0,h)}}
=\frac{\dsp\scal{\mathfrak{m}_{\alpha_{\bullet}}\phi_*}{\phi_*}_{L_2((0,\delta),h)}}
{\dsp\scal{\mathfrak{l}_{\alpha_{\bullet}}\phi_*}{\phi_*}_{L_2((0,\delta),h)}}%
=\Omega_*>0
$$
and, as explicit formulae above show,
\begin{equation}\label{eq:m_equal}
\begin{split}
\scal{\Mcal_\gamma\Psi}{\Psi}_{L_2(G^r_\gamma,h)}&=\scal{\Mcal_0\Psi}{\Psi}_{L_2(G^r_0,h)}
=\alpha_{\bullet}\int_{-r}^r\int_0^\delta \beta'(\eta) \erm^{-\beta(\eta)} |\phi_*(\eta)|^2\,\drm\eta\drm\xi\\
&= 2r\alpha_{\bullet}\int_0^\delta \beta'(\eta) \erm^{-\beta(\eta)} |\phi_*(\eta)|^2\,\drm\eta>0\,.
\end{split}
\end{equation}

We want to show that under conditions of Theorem~\ref{thm:main}
and with the choice of $\Psi$ as above, the inequality
\eqref{eq:main_ineq_mod} holds for any $r>R$.

In view of \eqref{eq:m_equal}, it is enough to show that
$$
D_\gamma:=\scal{\Lcal_\gamma\Psi}{\Psi}_{L_2(G^r_\gamma,h)}-\scal{\Lcal_0\Psi}{\Psi}_{L_2(G^r_0,h)}
$$
is \emph{negative} for $r>R$.

Explicit substitution gives, after taking into account the formula
$$
\int \frac{1}{p(\xi,\eta)^2}\,\frac{\dpar p(\xi,\eta)}{\dpar\xi}\,\drm\xi = 0
$$
(due to \eqref{eq:gam_compact}, with account of \eqref{eq:p}), the
following expression:
$$
D_\gamma = \int_{-r}^r\int_0^\delta \gamma(\xi)\eta\erm^{-\beta(\eta)}|\phi'_*(\eta)|^2\,\drm\eta\drm\xi
-\int_{-r}^r\int_0^\delta  \alpha_{\bullet}^2 \frac{\eta\gamma(\xi)}{1+\eta\gamma(\xi)}
\erm^{-\beta(\eta)}|\phi_*(\eta)|^2\,\drm\eta\drm\xi\,.
$$
This, in turn, can be re-written, using the obvious identity
$$
\frac{\eta\gamma(\xi)}{1+\eta\gamma(\xi)} = \eta\gamma(\xi) - \frac{\eta^2\gamma(\xi)^2}{1+\eta\gamma(\xi)}
$$
as
\begin{equation}\label{eq:D_gamma}
\begin{split}
D_\gamma &= \int_{-r}^r\gamma(\xi)\int_0^\delta
\eta\erm^{-\beta(\eta)}
\left(|\phi'_*(\eta)|^2-\alpha_{\bullet}^2|\phi_*(\eta)|^2\right)\,\drm\eta\drm\xi
\\
&+\alpha_{\bullet}^2\int_{-r}^r\int_0^\delta \frac{\eta^2\gamma(\xi)^2}{1+\eta\gamma(\xi)}
\erm^{-\beta(\eta)}|\phi_*(\eta)|^2\,\drm\eta\drm\xi\,.
\end{split}
\end{equation}

We shall deal with the two terms in \eqref{eq:D_gamma} separately.

The first one is more difficult. As \eqref{eq:alpha_crit} yields explicitly
$$
\alpha_{\bullet}^2=\frac{\dsp \int_0^\delta \erm^{-\beta(\eta)}|\phi'_*(\eta)|^2\,\drm\eta}%
{\dsp \int_0^\delta \erm^{-\beta(\eta)}|\phi_*(\eta)|^2\,\drm\eta}\,,
$$
we get
\begin{equation}\label{eq:I1}
\begin{split}
I_1&:=\int_0^\delta\eta\erm^{-\beta(\eta)}
\left(|\phi'_*(\eta)|^2-\alpha_{\bullet}^2|\phi_*(\eta)|^2\right)\,\drm\eta
\\
&=\frac{1}{\dsp \int_0^\delta \erm^{-\beta(\eta)}|\phi_*(\eta)|^2\,\drm\eta}
\times\left(\int_0^\delta\eta\erm^{-\beta(\eta)}|\phi'_*(\eta)|^2\,\drm\eta
\cdot
\int_0^\delta\erm^{-\beta(\eta)}|\phi_*(\eta)|^2\,\drm\eta
\right.\\
&\qquad\qquad\qquad\left.-
\int_0^\delta\erm^{-\beta(\eta)}|\phi'_*(\eta)|^2\,\drm\eta
\cdot
\int_0^\delta\eta\erm^{-\beta(\eta)}|\phi_*(\eta)|^2\,\drm\eta\right)
\,.
\end{split}
\end{equation}

We want to show that the term in brackets is negative under some reasonable assumptions:

\begin{lem}\label{lem:I1}
Assume that the conditions of theorem~\ref{thm:main} hold. Then
$I_1<0$.
\end{lem}

The proof of Lemma~\ref{lem:I1} uses the following
simple fact:\footnote{we are grateful to Daniel Elton for a useful suggestion helping to prove this Lemma.}

\begin{lem}\label{lem:Daniel}
Let $(a,b)\subset(0,+\infty)$ be a finite interval, and let a
function $g:(a,b)\to\Rbb$ be non-negative and non-increasing.
Then
$$
\left(\int_a^b x g(x)f(x)\,\drm x\right)\cdot\left(\int_a^b f(x)\,\drm x\right) - 
\left(\int_a^b g(x)f(x)\,\drm x\right)\cdot\left(\int_a^b xf(x)\,\drm x\right)\le 0
$$
for any non-negative function $f:(a,b)\to\Rbb$.
\end{lem}

\begin{proof}[Proof of Lemma~\ref{lem:Daniel}]
We have
\begin{align*}
&\left(\int_a^b x g(x)f(x)\,\drm x\right)\cdot\left(\int_a^b f(x)\,\drm x\right) - 
\left(\int_a^b g(x)f(x)\,\drm x\right)\cdot\left(\int_a^b xf(x)\,\drm x\right)\\
&\qquad=
\int_a^b\int_a^b x g(x)f(x) f(y)\,\drm x\drm y - \int_a^b\int_a^b g(x)f(x)yf(y)\,\drm x\drm y\\
&\qquad=
\int_a^b\int_a^y (x-y) f(x) f(y) g(x)\,\drm x\drm y + \int_a^b\int_y^b (x-y) f(x) f(y) g(x)\,\drm x\drm y
\end{align*}
Interchanging the variables $x$ and $y$ in the last integral, we obtain that the whole expression is equal to
$$
\int_a^b\int_a^y \underbrace{(x-y)}_{\text{non-positive}} \underbrace{f(x) f(y)(g(x)-g(y))}_{\text{non-negative}}\,\drm x\drm y\,,
$$
and is therefore non-positive.
\end{proof}

We can now proceed with evaluating $I_1$.

\begin{proof}[Proof of Lemma~\ref{lem:I1}] We act by doing a lot of integrations by part.
We shall also use \eqref{eq:phi_star_eq}.

We have (all integrals are over $[0,\delta]$ and with respect to $\eta$)
\begin{equation*}
\begin{split}
\int \eta\erm^{-\beta}|\phi'_*|^2 &= - \int \phi_*\cdot(\eta \erm^{-\beta} \phi'_*)'\\
&= -\int \phi_*\cdot(\erm^{-\beta} \phi'_*-\beta'\eta\erm^{-\beta}\phi'_*+\eta\erm^{-\beta} \phi''_*)\\
&=-\int \phi_*\cdot(\erm^{-\beta} \phi'_*+(\alpha_{\bullet}^2-\Lambda_*\alpha_{\bullet}\beta')\erm^{-\beta}\eta\phi_*)\,.
\end{split}
\end{equation*}
Further,
\begin{equation*}
-\int (\phi_*\erm^{-\beta})\phi'_*=\left(\int
\left(\phi'_*\erm^{-\beta}-\beta'\phi_*\erm^{-\beta}\right)\phi_*\right)-\erm^{-\beta(\delta)}\phi_*^2(\delta)\,,
\end{equation*}
thus producing
\begin{equation}\label{eq:int1}
\begin{split}
\int \eta\erm^{-\beta}|\phi'_*|^2 &= -\frac{1}{2} \int \beta'\phi_*^2\erm^{-\beta}
 \underbrace{-\frac{1}{2}\erm^{-\beta(\delta)}\phi_*^2(\delta)}_\text{negative constant}\\
& -\alpha_{\bullet}^2\int\eta\erm^{-\beta}|\phi_*|^2
 +\Lambda_*\alpha_{\bullet}\int\eta\beta'\erm^{-\beta}|\phi_*|^2
\end{split}
\end{equation}

Also,
\begin{equation}\label{eq:int2}
\begin{split}
\int \erm^{-\beta}|\phi'_*|^2 &= -\int \erm^{-\beta}\phi_*(-\beta'\phi'_*+\phi''_*)\\
&= -\int
\erm^{-\beta}\phi_*\left(-\beta'\phi_*'+\beta'\phi_*'+\alpha_{\bullet}^2\phi_*-\Lambda_*\alpha_{\bullet}\beta'\phi_*\right)\\
&=
-\int \erm^{-\beta}\phi_*^2\,(\alpha_{\bullet}^2-\Lambda_*\alpha_{\bullet}\beta')\,.
\end{split}
\end{equation}

Substituting \eqref{eq:int1} and \eqref{eq:int2} into \eqref{eq:I1}, and simplifying, we get
\begin{equation}\label{eq:I1_1}
\begin{split}
&I_1 \cdot \underbrace{\dsp\erm^{-\beta(\eta)}|\phi_*(\eta)|^2}_\text{positive integral}
=
\int \eta\erm^{-\beta}|\phi'_*|^2 \cdot \int \erm^{-\beta}|\phi_*|^2
- \int \erm^{-\beta}|\phi'_*|^2 \cdot \int \eta\erm^{-\beta}|\phi_*|^2\\
&\qquad=
\left(\underbrace{\left(-\frac{1}{2} \int \beta'\erm^{-\beta}|\phi_*|^2\right)}_{\text{negative as }\beta'>0}+
\underbrace{\left(-\frac{1}{2}\erm^{-\beta(\delta)}\phi_*^2(\delta)\right)}_{\text{negative constant}}\right)
\int \erm^{-\beta}|\phi_*|^2\\
&\qquad+\underbrace{\left(\Lambda_*\alpha_{\bullet}\right)}_{\text{$+$ve constant}}
\times
\underbrace{\left(\int \eta\beta'\erm^{-\beta}|\phi_*|^2\cdot \int \erm^{-\beta}|\phi_*|^2-
\int \beta'\erm^{-\beta}|\phi_*|^2\cdot \int \eta\erm^{-\beta}|\phi_*|^2\right)}_{\text{non-positive by Lemma~\ref{lem:Daniel} with }
g=\beta',\ f=\erm^{-\beta}|\phi_*|^2\text{ as }g'=\beta''\le 0}
\end{split}
\end{equation}
Thus $I_1<0$.
\end{proof}

Let us now return to \eqref{eq:D_gamma} and deal with the second term in the right-hand side. We have,
with account of \eqref{eq:gam_ext} and \eqref{eq:kappa_rest},
$$
\frac{\eta^2\gamma(\xi)^2}{1+\eta\gamma(\xi)}\le
\begin{cases}
\eta^2\gamma(\xi)^2\,,\qquad&\text{if }\gamma(\xi)\ge 0\\
\dsp\frac{1}{1-A}\eta^2\gamma(\xi)^2\,,\qquad&\text{if }\gamma(\xi)< 0
\end{cases}\qquad\le\max\{1,\dsp\frac{1}{1-A}\}\eta^2\gamma(\xi)^2\,,
$$
and so
$$
\alpha_{\bullet}^2\int_{-r}^r\int_0^\delta \frac{\eta^2\gamma(\xi)^2}{1+\eta\gamma(\xi)}
\erm^{-\beta(\eta)}|\phi_*(\eta)|^2\,\drm\eta\drm\xi\le
I_2\int_{-r}^r \gamma(\xi)^2\drm\xi\,,
$$
where
\begin{equation}\label{eq:I2}
I_2:=\max\{1,\dsp\frac{1}{1-A}\}\alpha_{\bullet}^2\int_0^\delta \eta^2\erm^{-\beta(\eta)}|\phi_*(\eta)|^2\,\drm\eta\ge 0\,.
\end{equation}

Thus, as $\gamma(\xi)$ vanishes for $|\xi|>R$, we have
$$
D_\gamma=I_1\int \gamma(\xi)\drm\xi+I_2\int \gamma(\xi)^2\drm\xi=
(-I_1)\left(C_\beta\int \gamma(\xi)^2\drm\xi-\int \gamma(\xi)\drm\xi\right)\,,
$$
where 
\begin{equation}\label{eq:C_beta}
C_\beta=\frac{I_2}{-I_1}=\frac{\dsp\max\{1,\dsp\frac{1}{1-A}\}\alpha_{\bullet}^2\int_0^\delta \eta^2\erm^{-\beta(\eta)}|\phi_*(\eta)|^2\,\drm\eta}
{\dsp\int_0^\delta\eta\erm^{-\beta(\eta)}
\left(|\phi'_*(\eta)|^2-\alpha_{\bullet}^2|\phi_*(\eta)|^2\right)\,\drm\eta}
\end{equation}
is a positive constant.

As soon as \eqref{eq:main_cond} holds, $D_\gamma$ is negative, and so \eqref{eq:main_ineq_mod} holds. This proves Theorem~\ref{thm:main}.

Finally, it is sufficient to note that \eqref{eq:cor_cond} implies $\dsp\int \gamma(\xi)^2\drm\xi<2Rc_{\beta,R}\int \gamma(\xi)\drm\xi$,
which proves Corollary~\ref{cor:main}.

\section{Conclusions}
It has been shown that a trapped mode is possible in the model presented here.
To increase confidence that such modes exist on real coasts further work is clearly
required to demonstrate that this mode is not an artifact of the modelling assumptions.
However these assumptions are the usual ones for the simple theory of CSWs and extensions
to include stratification and more realistic boundary conditions have not in general
contradicted them \cite{LeBloM78}.
The result here suggests that it would be of interest to compare low frequency
velocity records in the neighborhood of capes with those on nearby straight
coasts to determine whether there is indeed enhanced energy at the cape. Both
the above endeavors are being pursued.

\bibliographystyle{alpha}

\end{document}